\documentclass{article}
\usepackage{amsmath,amssymb,amsthm}
\usepackage{graphics}

\newtheorem{theorem}{Theorem}

\title {Circle actions with two fixed points}

\author {Oleg R. Musin\thanks{This research is partially supported by the NSF grant DMS-1400876 and the RFBR grant 15-01-99563.}}

\begin{document}
\date{}
\maketitle

\begin{abstract} We prove that if the circle group acts smooth and unitary on $2n$--dimensional stably complex manifold with two isolated fixed points and it is not bound equivariantly, then $n=1$ or $3$.   Our proof relies on the rigid Hirzebruch  genera. 
\end {abstract}

\section{Introduction}
We say that $M$ is a  {\it unitary ${\bf S}^1$--manifold} if it is a smooth
closed manifold with an effective circle action such that its tangent
bundle admits an ${\bf S}^1$--equivariant stable complex structure.

Let $M$ be a unitary $2n$--dimensional ${\bf S}^1$--manifold with $m$ isolated fixed points $p_1,\ldots,p_m$. Denote  weights in $p_i$ by $w_{i1},\ldots,w_{in}$ and its sign by $\varepsilon_i$, where $\varepsilon_i=\pm 1$. 

If $m=2$, then known the following sets of weights for unitary ${\bf S}^1$--manifolds:

$(Z)$: $w_{1i}=w_{2i}=a_i$, where $a_i\in{\Bbb Z}, \, a_i\ne0,$ for all $i=1,\ldots,n$, and  $\varepsilon_1=-\varepsilon_2$. 
A manifold  with these weights is cobordant to zero in the ring $U_*^{{\bf S}^1}$ of complex bordisms  with circle actions.  

($L_1$): $n=1$, $w_{11}=a, \, w_{21}=-a$, where $a>0$, and $\varepsilon_1=\varepsilon_2$. It is a linear action of ${\bf S}^1$ on ${\Bbb C}${P}$^1$:
 $$
 [z_0:z_1] \to  [e^{ia\varphi}z_0:z_1]. 
 $$

$(S_3)$: $n=3$, weights in $p_1$ and $p_2$ are $(a,b,-(a+b))$ and $(-a,-b,a+b)$ respectively, where $a$ and $b$ are positive integer, and $\varepsilon_1=\varepsilon_2$. This action admit ${\Bbb S}^6$ that is almost complex.  

Our main result is:

\begin{theorem}
Let $M$ be a unitary $2n$--dimensional ${\bf S}^1$--manifold with two isolated fixed points. Then the set of weights can only be $Z$, $L_1$ or $S_3$.
\end{theorem} 

Our proof relies on the fact that $T_{x,y}$ genus is rigid. We consider this method in the next section. 

 Kosniowski  \cite{Kos} conjectured that a unitary $2n$--dimensional ${\bf S}^1$--manifold that is not bound equivariantly has at least $\lceil{n/4}\rceil+1$ isolated fixed points.  Theorem 1 proves the conjecture for $n\le 4$. Note that this problem was solved only for ${\bf T}^n$--manifolds. Namely,  Zhi L\"u \cite{Lu} was proved that for torus manifolds there are at least $\lceil{n/2}\rceil+1$ isolated fixed points. 

It is not clear, can Kosniowski's conjecture be proved using only the rigidity of the $T_{x,y}$ genus? Perhaps,  explicit generators of $U_*^{{\bf S}^1}$ from our paper \cite{mus} can be useful for that.

\section{Rigid Hirzebruch genera}

Let $U_*$ be the complex bordism ring with coefficients in $R={\Bbb Q}$, ${\Bbb R}$, or ${\Bbb C}$. For a closed smooth stably complex manifold $X$, Hirzebruch \cite{Hir1} defined a multiplicative genus $h(X)$ by a homomorphism $h:U_*\otimes R\to R$.

Recall that according to Milnor and Novikov, two stably complex  manifolds are complex cobordant if and only if they have the same Chern numbers. Therefore, for any  multiplicative genus $h$ there exists a multiplicative sequence of polynomials $\{K_i(c_1,\ldots,c_i)\}$  such that
$
h(X)=K_n(c_1,\ldots,c_n),
$
where the $c_k$ are the Chern classes of $X$ and $n=\dim_{\Bbb C}(X)$ (see \cite{Hir1}).

Let $U_*^G$ be the ring of complex bordisms of manifolds with actions of a compact Lie group $G$. Then for any homomorphism $h:U_*\otimes R\to R$ we can define the {\it equivariant genus} $h^G$, i.e. a homomorphism
$$
h^G:U_*^G\otimes R\to K(BG)\otimes R
$$
(see details in \cite{krich1}).

A multiplicative genus $h$ is called {\it rigid} if for any connected compact group $G$ the equivariant genus $h^G(X)=h(X)$. For the complex case rigidity means that
$$
h^G:U_*^G\otimes R\to R\subset K(BG)\otimes R,
$$
i.e. $h^G([X,G])$ belongs to the ring of constants. It is well known (see \cite{AH, krich1}) that ${\bf S}^1$-rigidity implies $G$-rigidity, i.e. it is sufficient to prove rigidity only for the case $G={\bf S}^1$.

For $G={\bf S}^1$, the universal classifying space $BG$ is ${\Bbb C}${P}$^\infty$, and the ring $K(BG)\otimes R$ is isomorphic to the ring of formal power series  $R[[u]]$. Then for any ${\bf S}^1-$manifold $X$ and a Hirzebruch genus $h$ we have 
$h^{S^1}([X,{\bf S}^1])$ in $R[[u]]$. In the case of an ${\bf S}^1-$action on a unitary manifold $X^{2n}$ with isolated fixed points $p_1\ldots,p_m$ with weights $w_{i1},\ldots,w_{in}$, and signs $\varepsilon_i$, $i=1,\ldots,m$,  $\; h^{S^1}$ can be found explicitly (see \cite{BR,krich1}):
$$
h^{S^1}([X,{\bf S}^1])=S_h(\{w_{ij}\},u):=\sum\limits_{i=1}^m\varepsilon_i\prod_{j} {\frac{H(w_{ij}\,u)}{w_{ij}u}}, \eqno (1) 
$$
where $H$ is the characteristic series of $h$ (see \cite{Hir1,mus}). If $h$ is rigid, then from $(1)$ it follows that
$$
h(X)=S_h(\{w_{ij}\},u) \; \mbox{ for any } \; u. \eqno (2)
$$

Atiyah and based on the Atiyah-Singer index theorem  proved that  $T_y-$genus is rigid \cite{AH}.  Krichever \cite{krich1} gives a proof of rigidity of the $T_{x,y}-$genus using global analytic properties of $S_h(\{w_{ij}\},u)$. It is not hard to see that $(2)$ yields the
Atiyah-Hirzebruch formula \cite{AH,krich1, BR}  for a unitary ${\bf S}^1$-manifold $X$:
$$
T_{x,y}(X)=\sum\limits_{i=1}^m {\varepsilon_ix^{s_i^+}(-y)^{s_i^-}}, \eqno (3)
$$
where $s_i^+, s_i^-$ are numbers of positive and negative weights $\{w_{ij}\}$.

In \cite{krich1} is given the following formula for the $T_{x,y}$--genus: 
 $$
 H_{x,y}(u)=\frac{u(xe^{u(x+y)}+y)}{e^{u(x+y)}-1}.
 $$

Let $z:=e^{(x+y)u}$. Then (1), (2), and (3) imply
$$
\sum\limits_{i=1}^m\varepsilon_i\prod_{j} {\frac{xz^{w_{ij}}+y}{z^{w_{ij}}-1}}\equiv \sum\limits_{i=1}^m {\varepsilon_ix^{s_i^+}(-y)^{s_i^-}} \eqno (4) 
$$

 \section{Proof of the main theorem}
\begin{proof} Since we have only two fixed points the graph of weights (see \cite{mus80}) has two vertices and $n$ edges join these vertices. It implies that $|w_{1i}|=|w_{2i}|$ for all $i=1,\ldots,n.$ 
(Actually, it can be derived from (4), where $x=y=1$, see \cite[Theorem 1.1]{mus80}.) Moreover, if the set of weights is not $Z$, then $w_{1i}=-w_{2i}$. Otherwise, $S_{T_{x,y}}(\{w_{ij}\},u)$ cannot be constant. Therefore, for this case formula (4) can be written in the following form: 
$$
\frac{(xz^{a_1}+y)\cdots(xz^{a_k}+y)(x+yz^{b_1})\cdots(x+yz^{b_\ell})}
{(z^{a_1}-1)\cdots(z^{a_k}-1)(z^{b_1}-1)\cdots(z^{b_\ell}-1)} - 
$$
$$
\frac{(x+yz^{a_1})\cdots(x+yz^{a_k})(xz^{b_1}+y)\cdots(xz^{b_\ell}+y)}
{(z^{a_1}-1)\cdots(z^{a_k}-1)(z^{b_1}-1)\cdots(z^{b_\ell}-1)}=x^ky^\ell-x^\ell y^k, \eqno (5) 
$$
where $k+\ell=n$ and all $a_i$ and $b_j$ are positive integers.

So for $n=1$ we have $L_1$. If $n>1$, then $k>0$ and $\ell>0$. For $y=0$ (5) implies
$$
a_1+\ldots+a_k=b_1+\ldots+b_\ell \eqno (6) 
$$

Without loss of generality we may assume that $a:=a_1\ge a_i$ for all $i$ and $a_1>b_j$ for all $j$. (As we mentioned above, we cannot have the equality $a_1=b_j$.) 
Let $x=-z^a$ and $y=1$. Then (5) yields 
$$
\frac{(z^{2a}-1)\cdots(z^{a_k+a}-1)(z^{a}-z^{b_1})\cdots(z^{a}-z^{b_\ell})}
{(z^{a}-1)\cdots(z^{a_k}-1)(z^{b_1}-1)\cdots(z^{b_\ell}-1)}
=(-1)^\ell z^{ak}-(-1)^k z^{a\ell} \eqno (7) 
$$
Since $a>b_j$, we have $a\ell>b_1+\ldots+b_\ell$. Then (7) implies 
$$
ka=b_1+\ldots+b_\ell \eqno (8) 
$$
Therefore, from (6) we have $a_1+\ldots+a_k=ka$, $a_1=\ldots=a_k=a$, $k$ is odd, and $\ell$ is $even$. Then 
$$
\frac{(z^{a}+1)^k(z^{a-b_1}-1)\cdots(z^{a-b_\ell}-1)}
{(z^{b_1}-1)\cdots(z^{b_\ell}-1)}
= z^{a(\ell-k)}+1 \eqno (9) 
$$
Equation (9) yields $k=1$, $\ell=2$ and $a=b_1+b_2$. Thus, it is the case $S_3$. 
\end{proof}

\medskip

\medskip

\medskip








\medskip

\medskip

\medskip

\medskip

 O. R. Musin,  University of Texas Rio Grande Valley, School of Mathematical and Statistical Sciences, One West University Boulevard, Brownsville, TX, 78520.

 {\it E-mail address:} oleg.musin@utrgv.edu


\begin{thebibliography}{99}

\bibitem{AH}
M. F. Atiyah and  F. Hirzebruch, Spin manifolds and group actions, in Essays in Topology and Related Subjects, Springer-Verlag, Berlin, 1970, pp. 18-28.




\bibitem{BR}
V. M. Buchstaber and N. Ray,  The universal equivariant genus and Krichever's formula, {\it Russian Math. Surveys}, 62:1 (2007), 178 - 180.


\bibitem{Hir1}

F. Hirzebruch, Topological Methods in Algebraic Geometry, 3rd edition, Grundlehren der mathematischen Wissenschaften, no. 131, Springer, Berlin-Heidelberg 1966.

\bibitem{Kos}
C. Kosniowski, Some formulae and conjectures associated
with circle actions, in ``Topology Symposium, Siegen 1979'' (Proc.
Sympos., Univ. Siegen, Siegen, 1979), 331--339, Lecture
Notes in Math., 788, Springer, Berlin, 1980.


\bibitem{krich1}
I. M. Krichever, Formal groups and the Atiyah - Hirzebruch formula, {\it Izv. Akad. Nauk SSSR Ser. Mat.} {\bf 38} (1974), 1289-1304.



\bibitem{Lu}
Z. L\"u, Equivariant cobordism of unitary toric manifolds, in Conference V.  Buchstaber--70  ``Algebraic topology and Abelian functions'', 18--22 June 2013, Moscow, Abstracts, 53--54.

\bibitem{mus11}
O. R. Musin, On rigid Hirzebruch genera, {\it Mosc. Math. J.,} {\bf 11:1} (2011), 139--147


\bibitem{mus}
O. R. Musin, Generators of ${\bf S}^1$-bordisms,  {\it Math. USSR Sb.} {\bf 44} (1983),  325-334

\bibitem{mus80}
O.R. Musin, Circle action on homotopy complex projective spaces, {\it Mathematical Notices,}  {\bf 28:1} (1980), 533--540. 



\end{thebibliography}
\end{document}